\UseRawInputEncoding

\documentclass[10pt]{amsart}

\newtheorem{lemma}{Lemma}

\newtheorem{prop}{Proposition}

\usepackage{amsmath}
\usepackage[psamsfonts]{amssymb}
\usepackage[mathscr]{euscript}

\usepackage{color}

\newcommand{\cal}{\mathcal}

\title[On a natural map between twistor spaces \hfill] {On a natural map
between twistor spaces}

\author{Johann Davidov}
\thanks{The author is partially supported by the National Science Fund, Ministry of Education and Science of Bulgaria under contract KP-06-N52/3
.}

\address{Johann Davidov\\Institute of Mathematics and Informatics \\
Bulgarian Academy of Sciences\\ Acad. G. Bonchev St. Bl. 8\\ 1113 Sofia\\
Bulgaria.}

\begin{document}

\begin{abstract}

A diffeomorphism between the twistor spaces of two Riemannian
metrics on a smooth manifold preserving the fibres is defined based
on a well-known construction. It is shown that this bundle
isomorphim is a holomorphic map with respect to the
Atiyah-Hitchin-Singer, respectively Eells-Salamon, almost complex
structure if and only if the two metrics are conformal, respectively
homothetic. In these cases, the diffeomorphism is the identity map
and the result obtained provides an interpretation of the well-known
fact that the Atiyah-Hitchin-Singer almost complex structure of a
twistor space is invariant under conformal changes of the metric on
the base manifold, while that of Eells-Salamon is not invariant in
general. The more general problem of when an arbitrary bundle
isomorphism between twisor spaces is holomorphic is also considered.
Another problem discussed in the paper is when the diffeomorphism
mentioned above is a harmonic map with respect to natural families
of Riemannian metrics on the twistor spaces defined by means of the
two Riemannian metrics. It is proved that if the metrics are
conformal this happens if and only if they are homothetic.

\vspace{0,1cm} \noindent 2010 {\it Mathematics Subject
Classification} 53C28; 32L25, 58E20.

\vspace{0,1cm} \noindent {\it Key words: twistor spaces, almost
complex structures, holomorphic maps, harmonic maps}.

\end{abstract}

\maketitle

\section{Introduction}

The twistor theory has its origin in Mathematical Physics. Inspired
by the Penrose program \cite{Pe1,Pe2}, Atiyah, Hitchin and Singer
\cite{AHS} have developed this theory on oriented Riemannian
$4$-manifolds. It has been further generalized to any even
dimensional Riemannian manifold by L. B\'erard-Bergery and T. Ochiai
\cite{BBO},  and  I. Skornyakov \cite{Sk},   see also N. O'Brian and
J. Rawnsley \cite{OR}), and S. Salamon \cite{Sal}.

The twistor space of a Riemannian manifold $M$ admits two natural
almost complex structures, which have been introduced  by
Atiyah-Hitchin-Singer \cite{AHS} and Eells-Salamon \cite{ES} in the
case $dim\,M=4$. These structures have been studied by many authors
in various aspects and have important applications in geometry.

The motivation of this paper comes from the well-known fact that the
Atiyah-Hitchin-Singer almost complex structure  on the twistor space
of a Riemannian manifold is invariant under conformal changes of the
metric on the manifold. On the other hand, the Eells-Salamon almost
complex structure  is not invariant in general. It is a natural
question whether there are relations of these almost complex
structures on the twistor spaces corresponding to two arbitrary
Riemannian metrics on a manifold.

Recall that the twistor space of a Riemannian manifold $(M,g)$ is
the  bundle $\pi:{\mathcal Z}\to M$ whose fibre at a point $p\in M$
consists of all $g$-orthogonal complex structures $I_p:T_pM\to
T_pM$, $I_p^2=-Id$, on the tangent space of $M$ at $p$. The fibre of
the bundle ${\mathcal Z}$ can be identified with  the compact
Hermitian symmetric space $O(2m)/U(m)$, $2m=dim\,M$. Suppose that
$\widetilde g$ is another Riemannian metric on $M$, and
$\widetilde{\mathcal Z}$ is the twistor space of the Riemannian
manifold $(M,\widetilde g)$. It is well-known that given a metric
and an  arbitrary almost complex structure on a manifold, one can
construct an almost complex structure which is orthogonal with
respect to the given metric. It seems that it is not known who first
made this constriction (it can be considered as "folklore"). In the
present paper, it is used in order to define a diffeomorphism
$\Psi:{\mathcal Z}\to \widetilde{\mathcal Z}$ preserving the fibres
as follows. Let $C$ be the endomorphism of the tangent bundle $TM$
determined by the formula $g(CX,Y)=\widetilde g(X,Y)$, $X,Y\in TM$.
The endomorphism $C$ is symmetric and positive with respect to both
metrics $g$ and $\widetilde g$, and we denote its principal square
root by $Q$. Then the diffeomeorphism $\Psi$ is denied by
$\Psi(I)=Q^{-1}IQ$, $I\in{\mathcal Z}$.

In Sec. 2 of the paper, we collect some basic facts about twistor
spaces.

In the first part of Sec. 3, we discuss the problem of when the
diffeomeorphism $\Psi$ preserves the Atiyah-Hitchin-Singer or
Eells-Salamon almost complex structures on ${\mathcal Z}$ and
$\widetilde{\mathcal Z}$, i.e., when $\Psi$ is a (pseudo)
holomorphic map. It is not hard to see that $\Psi$ is holomorphic
with respect to the Atiyah-Hitchin-Singer almost complex structures
if and only if the metrics $g$ and $\widetilde g$  are conformally
equivalent. In this case, ${\mathcal Z}=\widetilde{\mathcal Z}$
(clearly, every $g$-orthogonal complex structure is $\widetilde
g$-orthogonal and vice versa) and $\Psi=Id$. Thus,  we have an
interpretation of the fact mentioned above that the
Atiyah-Hitchin-Singer almost complex structure is invariant under
conformal changes of the metric on the base manifold. We also
observe that $\Psi$ is holomorphic with respect to the Eells-Salamon
almost complex structures exactly when the metrics $g$ and
$\widetilde g$ are homothetic. The more general problem of when an
arbitrary bundle isomorphism between twisor spaces is holomorphic is
discussed in this section as well.

The twistor space ${\mathcal Z}$ of a Riemannian manifold $(M,g)$
admits a natural $1$-parameter family of Riemannian metrics $g_s$,
$s>0$, such that the projection map $\pi:({\mathcal Z},g_s)\to
(M,g)$ is a Riemannian submersion with totally geodesic fibres. The
Atiyah-Hitchin-Singer and Eells-Salamon almost complex structures
are orthogonal with respect to each metric $g_s$. Thus, we obtain
two families of almost Hermitian structures. These structures have
interesting geometric properties studied in numerous works.

In the second part of Sec. 3, we consider the problem of when the
map $\Psi:({\mathcal Z},g_s)\to (\widetilde{\mathcal Z},\widetilde
g_t)$ is harmonic. After computing the second fundamental form of
$\Psi$, we give a complete answer to this question in the most
interesting for us case of conformal  metrics $g$ and $\widetilde
g$. It is proved that, in this case, $\Psi=Id$ is harmonic if and
only if $g$ and $\widetilde g$ are homothetic.

\section{Basics of twistor spaces}

In this section, we collect some facts about twistor spaces that
will be used further on.

\subsection{The manifold of compatible linear complex structures}\label{CLCS}

Let $V$ be a real vector space of even dimension $n=2m$ endowed with
an Euclidean metric $g$. Denote by $F(V)$ the set of all complex
structures on $V$ compatible with the metric $g$, i.e.,
$g$-orthogonal. This set possesses the structure of an embedded
submanifold of the vector space $\mathfrak{so}(V)$ of skew-symmetric
endomorphisms of $(V,g)$. The tangent space of $F(V)$ at a point $J$
consists of all endomorphisms $Q\in \mathfrak{so}(V)$ anti-commuting
with $J$, and we have the decomposition
\begin{equation}\label{decom}
\mathfrak{so}(V)=T_JF(V)\oplus \{S\in\mathfrak{so}(V):~ SJ-JS=0\}.
\end{equation}
This decomposition is orthogonal with respect to the metric
$G(S,T)=-\displaystyle{\frac{1}{2}}\mathit{Trace}\, (S\circ T)$ of
$\mathfrak{so}(V)$.  Clearly the orthogonal projection ${\mathcal
V}_{J}\varphi$ of an endomorphism $\varphi\in \mathfrak{so}(V)$ to
the tangent space $T_JF(V)$ is given by
\begin{equation}\label{proj}
{\mathcal V}_J\varphi=\frac{1}{2}(\varphi+J\circ \varphi \circ J).
\end{equation}

 The smooth manifold $F(V)$ admits a natural
almost complex structure ${\mathcal J}$ defined by ${\mathcal
J}Q=J\circ Q$ for $Q\in T_JF(V)$. This structure is compatible with
the metric $G$. It is not hard to find an explicit formula for the
Levi-Civita connection $\nabla$ of $G$ and to see that $\nabla
{\mathcal J}=0$, so $(G,{\mathcal J})$ is a K\"ahler structure on
$F(V)$. In particular, ${\mathcal J}$ is an integrable structure.
Note also that the metric $G$ is Einstein with scalar curvature
$\frac{1}{2}m(m-1)^2$, see, e.g., \cite{D}.

The group $O(V)\cong O(2m)$ of orthogonal transformations of $V$
acts on $F(V)$ by conjugation and the isotropy subgroup at a fixed
complex structure $J_0$ is isomorphic to the unitary group $U(m)$.
Therefore $F(V)$ can be identified with the homogeneous space
$O(2m)/U(m)$. In particular, $dim\,F(V)=m^2-m$. Note also that the
manifold $F(V)$ has two connected components $F_{\pm}(V)$: if we fix
an orientation on $V$, these components consists of all complex
structures on $V$ compatible with the metric $g$ and inducing $\pm$
the orientation of $V$; each of them has the homogeneous
representation $SO(2m)/U(m)$.

The metric $g$ of $V$ induces a metric on $\Lambda^2V$ given by
\begin{equation}\label{g on Lambda2}
g(v_1\wedge v_2,v_3\wedge
v_4)=\frac{1}{2}[g(v_1,v_3)g(v_2,v_4)-g(v_1,v_4)g(v_2,v_3)].
\end{equation}
Then we have an isomorphisms $\mathfrak{so}(V)\cong \Lambda^2V$,
which sends $\varphi\in \mathfrak{so}(V)$ to the $2$-vector
$\varphi^{\wedge}$ determined by the identity
\begin{equation}\label{wedge-isom}
2g(\varphi^{\wedge},u\wedge v)=g(\varphi u,v),\quad u,v\in V.
\end{equation}
This isomorphism is an isometry with respect to the metric
$\frac{1}{2}G$ on $\mathfrak{so}(V)$ and the metric $g$ on
$\Lambda^2V$.

Suppose that $V$ is of dimension 4. Then the Hodge star operator
defines an endomorphism $\ast$ of $\Lambda^2V$ with $\ast^2=Id$.
Hence we have the orthogonal decomposition
$$
\Lambda^2 V=\Lambda^2_{+}V\oplus\Lambda^2_{-}V,
$$
where $\Lambda^2_{\pm}V$ are the subspaces of $\Lambda^2V$
corresponding to the $(\pm 1)$-eigenvalues of the operator $\ast$.
Let $(e_1,e_2,e_3,e_4)$ be an oriented orthonormal basis of $V$. Set
\begin{equation}\label{s-basis}
s_1^{\pm}=e_1\wedge e_2 \pm e_3\wedge e_4, \quad s_2^{\pm}=e_1\wedge
e_3\pm e_4\wedge e_2, \quad s_3^{\pm}=e_1\wedge e_4\pm e_2\wedge
e_3.
\end{equation}
Then $(s_1^{\pm},s_2^{\pm},s_3^{\pm})$ is an orthonormal basis of
$\Lambda^2_{\pm}V$. The orientation of $\Lambda^2_{\pm}V$ determined
by this basis does not depend on the choice of the basis
$(e_1,e_2,e_3,e_4)$ (see e.g. \cite{D17}).

The isomorphism $\varphi\to \varphi^{\wedge}$ identifies
$F_{\pm}(V)$ with the unit sphere $S(\Lambda^2_{\pm}V)$ of the
Euclidean vector space $(\Lambda^2_{\pm}V,g)$ (so, the factor $1/2$
in (\ref{g on Lambda2}) is chosen in order to have spheres with
radius $1$). Under this isomorphism, if $J\in F_{\pm}(V)$, the
tangent space $T_JF(V)=T_JF_{\pm}(V)$ is identified with the
orthogonal complement $({\mathbb R} J)^{\perp}$ of the space
${\mathbb R}J$ in $\Lambda^2_{\pm}V$.  Moreover, the restriction of
the complex structure ${\mathcal J}$ to $F_{\pm}(V)$ is identified
with the complex structure of the unit sphere in $\Lambda^2_{\pm}V$.

\subsection {The twistor space of a Riemannian manifold.} Let $(M,g)$
be a (connected) Riemannian manifold of even dimension.  The twistor
space of $(M,g)$ is the bundle $\pi:{\mathcal Z}\to M$ whose fibre
at a point $p\in M$ is the manifold $F(T_pM)$ of all complex
structures of the tangent space $T_pM$ compatible with  the metric
$g_p$ (i.e., $g_p$-orthogonal). According to this definition,  the
twistor spaces of two conformal metrics on a manifold coincide since
any complex structure on $T_pM$ compatible with one of these metrics
is compatible with the other one as well.

If $M$ is oriented, we can consider the bundles over $M$ whose fibre
at $p$ is the manifold of compatible complex structures yielding the
positive, respectively, the negative orientation of $T_pM$. The
total spaces of these  bundles are the connected components of the
manifold ${\mathcal Z}$ and are frequently called the positive and
the negative twistor spaces of $(M,g)$, respectively. They are
usually denoted by ${\mathcal Z}_{+}$ and ${\mathcal Z}_{-}$.
Clearly, changing the orientation interchanges the role of the
positive and negative twistor spaces.

If $M$ is oriented and of dimension 4, ${\mathcal Z}_{\pm}$ can be
identified with the unit sphere subbundle of $\Lambda^2_{\pm}TM$
with respect to the metric (\ref{g on Lambda2}), the latter vector
bundle being the subbundle of $\Lambda^2TM$ corresponding to the
$(\pm 1)$-eigenvalue of the Hodge star operator $\ast$. This
interpretation of a twistor space leads to different twistor spaces
of two conformally equivalent metrics $g$ and $\widetilde
g=e^{2f}g$; they are diffeomorphic by the map $\sigma\to
e^{-2f}\sigma$.

\subsection{The Atiyah-Hitchin-Singer  and Eells-Salamon almost complex structures on a twistor space.}
The manifold ${\mathcal Z}$ admits two almost complex structures
${\mathcal J}_1$ and ${\mathcal J}_2$ defined, respectively, by
Atiyah-Hitchin-Singer \cite{AHS} and Eells-Salamon \cite{ES} in the
case when the base manifold $M$ is of dimension 4.

The bundle ${\mathcal Z}$ can be considered as a subbundle of the
bundle $A(TM)$ of skew-symmetric endomorphisms of the tangent bundle
$TM$. Let $\nabla$ be a metric connection on $(M,g)$, and denote
 the induced connection on the bundle $Hom(TM,TM)$ again by $\nabla$.
The latter connection preserves the bundle $A(TM)$ since the
connection on $(M,g)$ is metric. The horizontal space of the vector
bundle $A(TM)$ with respect to $\nabla$ at a point $J\in{\mathcal
Z}$ is tangent to ${\mathcal Z}$. Hence we have the splitting
$T{\mathcal Z}={\mathcal H}\oplus {\mathcal V}$  of the tangent
bundle of ${\mathcal Z}$ into horizontal and  vertical subbundles.
The vertical subspace ${\mathcal V}_J$ of $T_J{\mathcal Z}$ is the
tangent space at $J$ to the fibre of the bundle ${\mathcal Z}$
through $J$. This fibre is the complex manifold $F(T_{\pi(J)}M)$ we
have discussed. Denote its complex structure by  ${\mathcal
J}_{{\mathcal V}_J}$. Then the almost complex structures ${\mathcal
J}_k$ are defined on the vertical space ${\mathcal V}_J$ by
$$
{\mathcal J}_1|{\mathcal V}_J={\mathcal J}_{{\mathcal V}_J},\quad
{\mathcal J}_2|{\mathcal V}_J=-{\mathcal J}_{{\mathcal V}_J}.
$$
For every $J\in {\mathcal Z}$, the horizontal subspace ${\mathcal
H}_J$ of $T_J{\mathcal  Z}$ is isomorphic to the tangent space
$T_{\pi(J)}M$ via the differential $\pi_{\ast\, J}$, and we set
$$
{\mathcal J}_k|{\mathcal H}_J=(\pi_{\ast}|{\mathcal H}_J)^{-1}\circ
J\circ (\pi_{\ast}|{\mathcal H}_J),\quad k=1,2.
$$

\noindent {\bf Notation}.  For $J\in {\mathcal Z}$ and $X\in
T_{\pi(J)}M$, we denote the horizontal lift $(\pi_{\ast}|{\mathcal
H}_J)^{-1}(X)$ of $X$ by $X^h_J$.

\hfill $\Box$

Note that if $X\in T_pM$ and $S$ is a section of ${\mathcal Z}$
considered as a section of the bundle $A(TM)$, then
$S_{\ast}(X)=X^h_{S(p)}+\nabla_{X}S$ where $\nabla_{X}S$ is a
skew-symmetric endomorphism of $T_pM$ anti-commuting with the
complex structure $S(p)$, i.e., a vertical tangent vector of
${\mathcal Z}$ at $S(p)$.

\subsection{Curvature decomposition} Let $R(X,Y,Z)=\nabla_{[X,Y]}Z-\nabla_{X}\nabla_{Y}Z
+\nabla_{Y}\nabla_{X}Z$ be the curvature tensor of the Levi-Civita
connection $\nabla$ of $(M,g)$. The curvature operator ${\mathcal
R}$ corresponding to the curvature tensor is the symmetric
endomorphism of $\Lambda ^{2}TM$ defined by
$$
g({\mathcal R}(X\wedge Y),Z\wedge U)=g(R(X,Y)Z,U),\quad X,Y,Z,U\in
TM.
$$
Denote the scalar curvature of the manifold $(M,g)$ by $\tau$, and
let $\rho:T_pM\to T_pM$ be its Ricci operator defined by
$g(\rho(X),Y)=Ricci(X,Y)$, $X,Y\in TM$. It is well-known that the
curvature operator ${\cal R}$ decomposes as (see, e.g., \cite[Sec. 1
G]{Besse}, \cite{ST69})
$$
{\mathcal R}=\frac{2\tau}{n(n-1)}Id+{\cal B}+{\cal W},
$$
where the operator ${\mathcal W}$ corresponds to the Weyl conformal
tensor $W$, while ${\cal B}$ corresponds to the traceless Ricci
tensor and is defined by
$$
{\mathcal B}(X\wedge Y)=\frac{2}{n-2}[X\wedge\rho(Y)+\rho(X)\wedge Y
-\frac{2\tau}{n}X\wedge Y], \quad X,Y\in T_pM.
$$
Note that this differs from \cite{Besse}  by a factor $\frac{1}{2}$
(the curvature operator therein is $\frac{1}{2}{\mathcal R}$)
because of the factor $1/2$ in our definition of the induced metric
on $\Lambda^2TM$.

The Riemannian manifold $(M,g)$ is Einstein exactly when ${\mathcal
B}=0$. It is locally conformally flat if and only if ${\mathcal
W}=0$. The manifold is of constant sectional curvature $c$ exactly
when ${\mathcal B}={\mathcal W}=0$ and $\tau=n(n-1)c$.

\smallskip

Suppose that $M$ is oriented and of dimension $4$. In this case, the
operator ${\mathcal B}$ sends $\Lambda^2_{\pm}TM$ into
$\Lambda^2_{\mp}TM$. Moreover,  the Weyl operator has an additional
decomposition ${\mathcal W}={\mathcal W}_{+}+{\mathcal W}_{-}$,
where $W_{\pm}$ is the restriction of ${\mathcal W}$ to
$\Lambda^2_{\pm}TM$ and zero on $\Lambda^2_{\mp}TM$, i.e.,
${\mathcal W}_{\pm}=\frac{1}{2}({\mathcal W}\pm\ast {\mathcal W})$
where $\ast$ is the Hodge star operator on $\Lambda^2TM$. The
manifold $(M,g)$ is called self-dual (anti-self-dual) if ${\mathcal
W}_{-}=0$ (resp. ${\mathcal W}_{+}=0$).

Note  that changing the orientation of $M$ interchanges the roles of
$\Lambda^2_{-}TM$ and $\Lambda^2_{+}TM$, hence the roles of
${\mathcal W}_{-}$ and ${\mathcal W}_{+}$.

Let $\widetilde g=e^{2f}g$ be a Riemannian metric conformal to $g$
and $\widetilde W$ its Weyl tensor. Then we have $\widetilde
g=e^{4f}g$ for the induced metrics on $\Lambda^2TM$. Hence
$\widetilde{\mathcal W}=e^{-2f}{\mathcal W}$ since $\widetilde
W(X,Y)Z=W(X,Y)Z$ for all $X,Y,Z\in TM$. Note also that the Hodge
star operator on $\Lambda^2TM$ does not change when the metric on
$M$ is replaced by a conformal one. Therefore the self-duality
(anti-self-duality) condition is invariant under conformal changes
of the metric. In connection with  this observation,  note that the
Atiyah-Hitchin-Singer almost complex structure ${\mathcal J}_1$ is
invariant under conformal changes of the metric of $M$, while, in
general, the Eells-Salamon one ${\mathcal J}_2$  is not.

\subsection {Integrability of the Atiyah-Hitchin-Singer  and
Eells-Salamon almost complex structures.} Let $\nabla$ be the
Levi-Civita connection of a Riemannian manifold $(M,g)$. It is a
result of Eells and Salamon \cite{ES, Sal} that the almost complex
structure ${\mathcal J}_2$  is never integrable. It has been proved
by Atiyah, Hitchin, and Singer \cite{AHS} that if the base manifold
$(M,g)$ is oriented and of dimension $4$, the almost complex
structure ${\mathcal J}_1$ on the negative twistor space is
integrable if and only if $(M,g)$ is self-dual. If $dim\,M\geq 6$,
${\mathcal J}_1$ is integrable exactly when the base manifold
$(M,g)$ is conformally flat, see,  e.g., \cite{Raw,Sal}. The
integrability conditions for the Atiyah-Hitchin-Singer  and
Eells-Salamon almost complex structures defined by means of an
arbitrary metric connection on $(M,g)$ have been discussed in
\cite{CGR}.

\smallskip

\subsection{A $1$-parameter family of Riemannian metrics on a
twistor space} The decomposition $T{\mathcal Z}={\mathcal H}\oplus
{\mathcal V}$
 allows one to define a natural $1$-parameter family of Riemannian
metrics $g_s$, $s>0$, on the manifold ${\mathcal Z}$ in the
following (standard) way.  The metric $ g_s$ on ${\mathcal H}_J$ is
the lift of the metric $g$ of $T_{\pi(J)}M$, $ g_s|{\mathcal
H}_J=\pi^{\ast} g$. On the vertical space ${\mathcal V}_J$, the
metric $g_s$ is defined as $s$ times the metric $G$ of the fibre
through $J$. Finally, the horizontal space ${\mathcal H}_J$ and the
vertical space ${\mathcal V}_J$ are declared to be orthogonal. Thus,
for $X,Y\in T_{\pi(J)}M$ and $V,W\in{\mathcal V}_J$,
\begin{equation}\label{ht}
g_s(X^h_J+V,Y^h_J+W)=g(X,Y)+sG(V,W).
\end{equation}

The almost complex structures ${\mathcal J}_1$ and ${\mathcal J}_2$
are both compatible with each metric $g_s$. Thus, for a fixed $s$,
we have two almost-Hermitian structures on the twistor space which
have been  studied by many authors.

\smallskip

For any $s>0$, the projection map $\pi:({\mathcal Z},g_s)\to (M,g)$
is a Riemannian submersion with totally geodesic fibres by the Vilms
theorem \cite{V}, \cite[Sec. 9 F]{Besse}.

Denote by $R$ the curvature of the connection on the vector  bundle
$A(TM)$ induced by the Levi-Civita connection $\nabla$ of $(M,g)$.

\begin{lemma}\label{R[a,b]} {\rm (\cite{D})} For every $a,b\in A(T_pM)$ and $X,Y\in T_pM$, we
have
$$
G(R(X,Y)a,b)=g(R([a,b]^{\wedge})X,Y).
$$
\end{lemma}

\begin{lemma} \rm (\cite{D,DM91})\label{Levi-Civita}
Let $D=D_s$ be the Levi-Civita connection of the metric $g_s$ on
${\mathcal Z}$. Let $J\in {\mathcal Z}$. If $X,Y$ are vector fields
on $M$ in a neighbourhood of the point $p=\pi(J)$ and $V$ is a
vertical vector field on ${\mathcal Z}$ near $J$, then
\begin{equation}\label{D-hh}
(D_{X^h}Y^h)_{J}=(\nabla_{X}Y)^h_{J}+\frac{1}{2}R_{p}(X\wedge Y)J
\end{equation}
\begin{equation}\label{D-vh}
(D_{V}X^h)_{J}={\mathcal H}(D_{X^h}V)_{J}=-s(R_{p}((J\circ
V_J)^{\wedge})X)^h_J \ ,
\end{equation}
where ${\mathcal H}$ means "the horizontal component".
\end{lemma}
\begin{proof} Identity (\ref{D-hh}) follows from the Koszul formula
for the Levi-Civita connection and the identity
$[X^h,Y^h]_J=[X,Y]^h_J+R_p(X,Y)J$, where $R_{p}(X,Y)J$ is a vertical
tangent vector of ${\mathcal Z}$ at $J$.

Let $W$ be a vertical vector field on ${\mathcal Z}$. Then
$$
g_s(D_{V}X^h,W)=-g_s(X^h,D_{V}W)=0 ,
$$
since the fibres are totally geodesic submanifolds, so $D_{V}W$ is a
vertical vector field. Therefore, the vector field $D_{V}X^h$ is
horizontal. Moreover, $[V,X^h]$ is a vertical vector field, hence
$D_{V}X^h={\mathcal H}D_{X^h}V$. Then at the point $J$, by
Lemma~\ref{R[a,b]},
$$
\begin{array}{c}
g_s(D_{V}X^h,Y^h)=g_s(D_{X^h}V,Y^h)=-g_s(V,D_{X^h}Y^h)=-\displaystyle{\frac{s}{2}}
G(V,R(X,Y)J)\\[6pt]
=-\displaystyle{\frac{s}{2}}g(R([J,V]^{\wedge})X,Y)=-s
g_s((R((J\circ V_J)^{\wedge})X)^h_J,Y^h_J).
\end{array}
$$
\end{proof}

\section{A natural map between twistor spaces}

Let $g$ and $\widetilde g$ be two Riemannian metrics on a manifold
$M$.

Define an endomorphism $C$ of the tangent bundle $TM$ setting
\begin{equation}\label{C}
g(C(X),Y)=\widetilde g(X,Y),\quad X,Y\in TM.
\end{equation}
For $p\in M$, $C(p)$ is a $g$-symmetric positive endomorphism of the
tangent space $T_pM$, and  we denote by $Q(p)$ the principal square
root of $C(p)$, i.e., the unique $g$-symmetric positive endomorphism
with $Q(p)^2=C(p)$, see, e.g., \cite[Theorem 7.2.6]{HJ}. It is well
known that $Q(p)$ depends smoothly on $p$ since $C(p)$ does so. For
reader's convenience, we give a proof of this fact. Note first that
by \cite[Theorem 1]{W65} for any symmetric positive definite matrix
$A$, the matrix equation $\exp X=A$ has a solution, denoted by $\ln
A$, with the integral representation
$$
\ln A=[A-I]\int_{0}^{1}[(1-\lambda)I+\lambda A]^{-1}d\lambda
$$
where $I$ stands for the identity matrix. This representation
implies that the matrix $\ln A$ is symmetric since $A$ is symmetric.
Then the matrix $B=\exp(1/2 \ln A)$ is symmetric, positive definite,
and such that $B^2=A$. Hence $B$ is the principal square root of
$A$. Furthermore, it follows from the integral representation of
$\ln A$ that the elements of  $\ln A$ are analytic functions of the
elements of $A$. Hence the elements of $B$ also possess this
property. Therefore if the elements of $A$ are smooth functions, so
are the elements of the principal square root. Now, take an
orthonormal frame of vector fields on $M$ and apply the latter
remark to the matrix representing the endomorphism $C$ to conclude
that $Q(p)$ depends smoothly on $p$.

Note that the endomorphisms $C(p)$ and $Q(p)$ are $\widetilde
g$-symmetric and positive since
$$
\begin{array}{c}
\widetilde g(C(X),Y)=g(CC(X),Y)=g(C(X),C(Y))=g(X,CC(Y))=\widetilde g(X,C(Y)), \\[6pt]

\widetilde g(Q(X),Y)=g(CQ(X),Y)=g(Q^3(X),Y)=g(QQ(X),Q(Y))\\[6pt]
\hspace {5cm}=g(X,CQ(Y))=\widetilde g(X,Q(Y)),
\end{array}
$$
$X,Y\in T_pM$, where, as usual, juxtaposition of endomorphisms means
their composition. Then clearly
\begin{equation}\label{g-h}
g(X,Y)=\widetilde g(Q^{-1}(X),Q^{-1}(Y)),\quad X,Y\in TM.
\end{equation}
Let $G(A,B)=\frac{1}{2}Trace\{X\to g(A(X),B(X))\}$  be the metric on
the vector bundle $Hom(TM,TM)$ of endomorphisms of $TM$ induced by
the metric $g$. Clearly, its restriction to the vertical spaces of
${\mathcal Z}$ coincides with the metric $G$ intoduced in
Sec.~\ref{CLCS}, hence the same notation. Similarly, denote by
$\widetilde G$ the metric on $Hom(TM,TM)$ induced by the metric
$\widetilde g$. Then, for $A,B\in Hom(TM,TM)$,
\begin{equation}\label{G-H}
G(A,B)=\widetilde G(Q^{-1}AQ,Q^{-1}BQ).
\end{equation}
Thus, the bundle map $\Psi:Hom(TM,TM)\to Hom(TM,TM)$ defined by
$\Psi(A)=Q^{-1}AQ$ is an isomorphism with respect to the metrics $G$
and $\widetilde G$. Moreover, an endomorphism $A$ of a tangent space
$T_pM$ is $g$-skew-symmetric (symmetric) if and only if $\Psi(A)$ is
$\widetilde g$-skew-symmetric (symmetric).

\smallskip
\noindent {\bf Notation}. The value $S(p)$ at a point $p$ of a
section $S$ of a bundle will also be denoted by $S_p$.

If $A$ and $B$ are sections of the bundle $Hom(TM,TM)$, the section
$M\ni p\to A_p\circ B_p$ of $Hom(TM,TM)$ will be denoted by $AB$.

As usual, the notations $\Psi\circ A$ and $\Psi(A)$ are used
interchangeably. \hfill $\Box$

\smallskip

Clearly, $\Psi(AB)=\Psi(A)\Psi(B)$ for any sections $A,B$ of
$Hom(TM,TM)$.

\smallskip

If $I$ is a complex structure on a tangent space $T_pM$ compatible
with the metric $g$, the endomorphism $\Psi(I)$ of $T_pM$ is a
complex structure on $T_pM$ compatible with the metric $\tilde g$
since
$$
\begin{array}{l}
{\widetilde g}(Q^{-1}IQ(X),Q^{-1}IQ(Y)) \\[6pt]
=g(CQ^{-1}IQ(X),Q^{-1}IQ(Y))=g(QIQ(X),Q^{-1}IQ(Y))\\[6pt]
=g(IQ(X),QQ^{-1}IQ(Y))=g(IQ(X),IQ(Y))=g(Q(X),Q(Y))=g(Q^2(X),Y)\\[6pt]
=\widetilde g (X,Y),
\end{array}
$$
where $Q$ stands for the endomorphism $Q(p)$ of $T_pM$.

\smallskip

The construction of the map $\Psi$  implies the well-known fact that
if $M$ is an almost complex manifold and $\widetilde g$ is any
metric on $M$, then the smooth manifold $M$ admits an almost complex
structure compatible with the metric $\widetilde g$. For, if $I$ is
a given almost complex structure on $M$, we can take any metric $g$
compatible with $I$, for example, $g(X,Y)=\widetilde
g(X,Y)+\widetilde g(IX,IY)$. Then $\Psi(I)$ is an almost complex
structure compatible with $\widetilde g$.

\smallskip

Suppose that the manifold $M$ is oriented. It is clear that if a
$g$-compatible complex structure $I$ on a tangent space $T_pM$
induces the orientation of $M$, the $\widetilde g$-compatible
complex structure $\Psi(I)$ also induces the orientation of $M$ (and
vice versa). Hence if ${\mathcal Z}$ and $\widetilde{\mathcal Z}$
denote the positive (or negative) twistor spaces of the Riemannian
manifolds $(M,g)$ and $(M,\widetilde g)$, we have a bundle
isomorphism $\Psi:{\mathcal Z}\to \widetilde{\mathcal Z}$.

\smallskip

\noindent {\bf Notation}. In what follows, ${\mathcal Z}$ will stand
for the positive twistor space of a Riemannian manifold.

\medskip

Let $\nabla$ and $\widetilde\nabla$ be the Levi-Civita connections
of the metrics $g$ and $\widetilde g$, respectively. Denote by
${\mathcal H}$ the horizontal subbundle of the tangent bundle
$T{\mathcal Z}$ defined by means of the connection  $\nabla$.
Similarly, $\widetilde{\mathcal H}$ will denote the horizontal
subbundle of $T\widetilde{\mathcal Z}$ defined via
$\widetilde\nabla$. If $X\in TM$, the horizontal lifts of $X$ to
${\mathcal H}$ and $\widetilde{\mathcal H}$ will be denoted by $X^h$
and $X^{\widetilde h}$, respectively. \hfill $\Box$

\smallskip

Recall that a smooth map $\varphi:(N,J)\to (\widetilde N,\widetilde
J)$ between almost complex manifolds is called holomorphic (or
pseudo-holomorphic) if $\widetilde J\circ
\varphi_{\ast}=\varphi_{\ast}\circ J$ (a generalization of the
Cauchy-Riemann equations); it is anti-holomorphic if $\widetilde
J\circ \varphi_{\ast}=-\varphi_{\ast}\circ J$.

\begin{prop}\label{J1}
Suppose the twistor spaces ${\mathcal Z}$ and $\widetilde{\mathcal
Z}$ are endowed with their Atiyah-Hitchin-Singer almost complex
structures  ${\mathcal J}_1$ and $\widetilde{\mathcal
J}_1$,respectively.

\smallskip

\noindent (i) The map $\Psi:({\mathcal Z},{\mathcal J}_1)\to
(\widetilde{\mathcal Z},\widetilde{\mathcal J}_1)$ is holomorphic if
and only if $C=e^{2f}Id$ for a smooth function $f$ on $M$, i.e., the
metrics $g$ and $\widetilde g$ are conformal with $\widetilde
g=e^{2f}g$. In this case ${\mathcal Z}=\widetilde {\mathcal Z}$ and
$\Psi=Id$.

\smallskip

\noindent (ii) The map $\Psi:({\mathcal Z},{\mathcal J}_1)\to
(\widetilde{\mathcal Z},\widetilde{\mathcal J}_1)$ is never
anti-holomorphic
\end{prop}

\begin{proof}
$(i)$ The map $\Psi$ sends a fibre of ${\mathcal Z}$ into a fibre of
$\widetilde{\mathcal Z}$, hence $\Psi_{\ast}$ maps a vertical vector
to a vertical vector. Let $I\in {\mathcal Z}$, and set $p=\pi(I)$
where $\pi$ is the projection map of the bundle ${\mathcal Z}\to M$.
Take a (local) section $S$ of ${\mathcal Z}$ such $S(p)=I$ and
$\nabla S|_p=0$. Then, if $X\in T_{p}M$,
\begin{equation}\label{diff on hor}
\Psi_{\ast}(X^h_I)=(\Psi\circ S)_{\ast}X=X^{\widetilde
h}_{\Psi(I)}+\widetilde\nabla_{X}(Q^{-1}SQ).
\end{equation}
Next, if $U\in T_I{\mathcal Z}$ is a vertical vector,
\begin{equation}\label{diff on ver}
\Psi_{\ast}(U)=Q^{-1}_p\circ U\circ Q_p.
\end{equation}
It follows from (\ref{diff on hor}) and (\ref{diff on ver}) that the
map $\Psi:{\mathcal Z}\to \widetilde{\mathcal Z}$ is holomorphic
with respect to the Atiyah-Hitchin-Singer almost complex structures
on ${\mathcal Z}$ and $\widetilde{\mathcal Z}$ if and only if for
every $I\in{\mathcal Z}$, $X\in T_{\pi(I)}M$, and $U\in{\mathcal
V}_I$,
\begin{equation}\label{cond-hol}
\begin{array}{c}
I(X)=\Psi(I)(X),\quad \widetilde\nabla_{I(X)}(Q^{-1} S Q)=\Psi(I)\circ\widetilde\nabla_{X}(Q^{-1} S Q),\\[6pt]
Q^{-1}_p\circ (I\circ U)\circ Q_p=\Psi(I)\circ (Q^{-1}_p\circ U\circ
Q_p),\quad p=\pi(I).
\end{array}
\end{equation}
The third identity of (\ref{cond-hol}) is trivially satisfied. The
first one is equivalent to $Q_{\pi(I)}\circ I=I\circ Q_{\pi(I)}$ for
every $I\in {\mathcal Z}$. Take an arbitrary oriented
$g$-orthonormal basis $e_1,...,e_{2m}$, $2m=n$, of a tangent space
$T_pM$. Let $I$ be the complex structure on $T_pM$ for which
$Ie_{2k-1}=e_{2k}$, $k=1,...,m$. Then $I\in{\mathcal Z}$, and the
identity $Q_{\pi(I)}\circ I=I\circ Q_{\pi(I)}$ gives
$g(Qe_1,e_2)=-g(Qe_2,e_1)$ and $g(Qe_1,e_1)=Q(e_2,e_2)$. The
endomorphism $Q$ is symmetric with respect to $g$, hence
$g(Qe_1,e_2)=0$. It follows that $g(Qe_i,e_j)=0$ and
$g(Qe_i,e_i)=g(Qe_j,e_j)$ for $i\neq j$, $i,j=1,...,n$. Therefore
$Q=e^{f}Id$ for a smooth function $f$ on $M$. Then the identity
$\widetilde\nabla_{I(X)}(Q^{-1} S
Q)=\Psi(I)\big(\widetilde\nabla_{X}(Q^{-1} S Q)\big)$ reduces to
\begin{equation}\label{cond-DS}
\widetilde\nabla_{I(X)}S=I\circ\widetilde\nabla_{X}S.
\end{equation}
The Levi-Civita connections of metrics $g$ and $\widetilde
g=e^{2f}g$ are related by
$$
\widetilde\nabla_{X}Y=\nabla_{X}Y +X(f)Y + Y(f)X-g(X,Y)\nabla f
$$
where $\nabla f$ is the gradient of $f$ defined by $g(\nabla
f,Z)=Z(f)$ (the dual to $df$ via $g$). Hence, for $X,Y\in T_pM$,
$$
(\widetilde\nabla_{X}S)(Y)=(\nabla_{X}S)(Y)+(IY)(f)X-Y(f)IX-g(X,IY)\nabla
f+g(X,Y)I(\nabla f).
$$
The latter formula implies (\ref{cond-DS}). This proves the first
claim.

\smallskip

\noindent $(ii)$ If $\Psi$ is anti-holomorphic, then for every $I\in
{\mathcal Z}$ and $U\in {\mathcal V}_I$,
\begin{equation}\label{not}
 Q^{-1}_p\circ (I\circ
U)\circ Q_p=-\Psi(I)\circ (Q^{-1}_p\circ U\circ Q_p),\quad p=\pi(I).
\end{equation}
This implies $IU=0$ for every $I$ and $U$, a contradiction.

\end{proof}

\begin{prop}\label{J_2} (i) The map $\Psi:({\mathcal Z},{\mathcal J}_2)\to
(\widetilde{\mathcal Z},\widetilde{\mathcal J}_2)$ is holomorphic if
and only if $\widetilde g=e^{2f}g$ for a constant $f$, i.e.,  the
metrics $g$ and $\widetilde g$ are homothetic (and $\Psi=Id$).

\smallskip

\noindent (ii) The map $\Psi:({\mathcal Z},{\mathcal J}_2)\to
(\widetilde{\mathcal Z},\widetilde{\mathcal J}_2)$ is never
anti-holomorphic.
\end{prop}

\begin{proof} We use the notation in the preceding proof.

If $\Psi$ is holomorphic or anti-holomorphic, then $I(X)=\Psi(I)(X)$
for every $I\in {\mathcal Z}$ and $X\in T_{\pi(I)}M$. As we have
noted, the latter identity implies $\widetilde g=e^{2f}g$ for a
smooth function $f$.

\smallskip

\noindent $(i)$ If $\Psi$ is holomorphic,
$\widetilde\nabla_{I(X)}(Q^{-1} S
Q)=-\Psi(I)\big(\widetilde\nabla_{X}(Q^{-1} S Q)\big)$. Hence
$\widetilde\nabla_{I(X)}S=-I\circ\widetilde\nabla_{X}S$. This
identity is equivalent to the identity
$$
g(IY,\nabla f)IX+g(Y,\nabla f)X-g(X,Y)\nabla f-g(X,JY)I(\nabla f)=0.
$$
Let $X\in T_pM$ be any vector orthogonal to $(\nabla f)(p)$ and
$I((\nabla f)(p))$. Then the above identity with $Y=X$ gives
$||X||^2\nabla f=0$. Hence $\nabla f=0$.

Conversely, it is easy to check that if $g$ and $\widetilde g$ are
homothetic, the map $\Psi:({\mathcal Z},{\mathcal J}_2)\to
(\widetilde{\mathcal Z},\widetilde{\mathcal J}_2)$ is holomorphic.

\smallskip

\noindent $(ii)$ The assumption that $\Psi$ is anti-holomorphic
implies identity (\ref{not})), which leads to a  contradiction.
\end{proof}

Similarly, we have the following.
\begin{prop}\label{J1-J_2}
The map $\Psi:{\mathcal Z}\to \widetilde{\mathcal Z}$ is neither
holomorphic, nor anti-holomorphic if one of the twistor spaces
${\mathcal Z}$ and $\widetilde{\mathcal Z}$ is endowed with  the
Atiyah-Hitchin-Singer almost complex structures and the other one
with the Eells-Salamon structure.
\end{prop}

Propositions~\ref{J1} - \ref{J1-J_2} can be placed in a more general
setting as follows.

\begin{prop}
Let ${\mathcal Z}$ and $\widetilde{\mathcal Z}$ be the twistor
spaces of Riemannian manifolds $(M, g)$ and $(\widetilde M,
\widetilde g)$. Let $\Phi:{\mathcal Z}\to\widetilde{\mathcal Z}$ be
a bundle isomorphism and  $\varphi:M\to \widetilde M$  the induced
diffeomorphism.

\smallskip

\noindent (i) The map $\Phi$ is holomorphic with respect to
${\mathcal J}_1$ and $\widetilde{\mathcal J}_1$ if and only if
$\Phi(I)=\varphi_{\ast}\circ I\circ\varphi_{\ast}^{-1}$. The latter
identity implies that the map $\varphi$ is conformal.

\smallskip

\noindent (ii) The map $\Phi$ is holomorphic with respect to
${\mathcal J}_2$ and $\widetilde{\mathcal J}_2$ if and only if
$\Phi(I)=\varphi_{\ast}\circ I\circ\varphi_{\ast}^{-1}$ and
$\varphi$ is a homothety.

\end{prop}

\begin{proof}
Let $I\in{\mathcal Z}$ and $p=\pi(I)$. Take a section $S$ of
${\mathcal Z}$ such that $S(p)=I$ and $\nabla S|_p=0$. Then, for
every $X\in T_pM$,
\begin{equation}\label{Phi-ast-hor}
\Phi_{\ast}X^h_I=(\Phi\circ S)_{\ast}(X)=(\Phi\circ
S\circ\varphi^{-1})_{\ast}(\varphi_{\ast}X)=(\varphi_{\ast}X)^{\widetilde
h}_{\Phi(I)}+\widetilde\nabla_{\varphi_{\ast}X}(\Phi\circ
S\circ\varphi^{-1}).
\end{equation}
It follows that the map $\Phi:({\mathcal Z},{\mathcal J}_k)\to
(\widetilde{\mathcal Z},\widetilde{\mathcal J}_k)$, $k=1$ or $2$, is
holomorphic if and only if for every $I\in{\mathcal Z}$ and $X\in
T_{\pi(I)}M$, and every $U\in{\mathcal V}_I$,
\begin{equation}\label{Jk-hol}
\begin{array}{c}
\Phi(I)(\varphi_{\ast}X)=\varphi_{\ast}(IX),\quad
(-1)^{k+1}\Phi(I)\circ \widetilde\nabla_{\varphi_{\ast}X}(\Phi\circ
S\circ\varphi^{-1})=\widetilde\nabla_{\varphi_{\ast}(IX)}(\Phi\circ
S\circ\varphi^{-1})\\[6pt]
\Phi(I)\circ \Phi_{\ast}(U)=\Phi_{\ast}(I\circ U)
\end{array}
\end{equation}
where $S$ is a local section of ${\mathcal Z}$ such that
$S(\pi(I))=I$ and $\nabla S|_{\pi(I)}=0$. Obviously, the first of
these identity is equivalent to
\begin{equation}\label{Phi}
\Phi(I)=\varphi_{\ast}\circ I\circ \varphi_{\ast}^{-1}.
\end{equation}
The latter identity implies $\Phi_{\ast}(U)=\varphi_{\ast}\circ
U\circ \varphi_{\ast}^{-1}$, and it is clear that the third identity
of (\ref{Jk-hol}) is satisfied. In order to show that, as a
consequence of (\ref{Phi}), the map $\varphi$ is conformal, let $q$
be a point of $M$, take a unit tangent vector $E\in T_qM$, and
consider an orthonormal basis $E_1,\cdots,E_n$ of $T_qM$ with
$E_1=E$. For any two different indexes $k\neq l$, there is a complex
structure $I$ of $T_qM$ such that $I\in {\mathcal Z}$ and
$IE_k=E_l$. The complex structure $\Phi(I)$ on
$T_{\varphi(q)}\widetilde M$ is compatible with the metric
$\widetilde g$, hence by (\ref{Phi})
$$
\widetilde g(\varphi_{\ast}E_k,\varphi_{\ast}E_l)=\widetilde
g(\Phi(I)\varphi_{\ast}E_k,\Phi(I)\varphi_{\ast}E_l)=\widetilde
g(\varphi_{\ast}(IE_k),\varphi_{\ast}(IE_l))=-\widetilde
g(\varphi_{\ast}E_l,\varphi_{\ast}E_k).
$$
Hence $\widetilde g(\varphi_{\ast}E_k,\varphi_{\ast}E_l)=0$.  Also,
taking $I\in{\mathcal Z}$ with $IE_1=E_k$, we see that $\widetilde
g(\varphi_{\ast}E_1,\varphi_{\ast}E_1)=\widetilde
g(\varphi_{\ast}E_k,\varphi_{\ast}E_k)$. Then, for
$X=\sum\limits_{i=1}^n\lambda_i E_i$ and $Y=\sum\limits_{i=1}^n\mu_i
E_i$,
$$
\widetilde
g(\varphi_{\ast}X,\varphi_{\ast}Y)=\sum\limits_{i=1}^n\lambda_i\mu_i\widetilde
g(\varphi_{\ast}E_1,\varphi_{\ast}E_1)=g(X,Y)\widetilde
g(\varphi_{\ast}E,\varphi_{\ast} E).
$$
It follows from this identity that the positive number $\widetilde
g(\varphi_{\ast}E,\varphi_{\ast} E)$ does not depend on the choice
of $E$. Moreover, setting $f(q)=\frac{1}{2}\ln \widetilde
g(\varphi_{\ast}E,\varphi_{\ast} E)$, we get a smooth function on
$M$ such that $\varphi^{\ast}\widetilde g=e^{2f}g$.

The Levi-Civita connection $\nabla^{\ast}$ of the metric
$\varphi^{\ast}\widetilde g$ coincides with the connection obtained
by transferring the Levi-Civita connection $\widetilde\nabla$ of
$\widetilde g$ via the diffeomorphism $\varphi$. Hence, if $X\in
T_qM$ and $Y$ is a vector field in a neighbourhood of $q$,
$$
\varphi_{\ast}(\nabla^{\ast}_{X}Y)=\widetilde\nabla_{\varphi_{\ast}X}(\varphi_{\ast}\circ
Y\circ\varphi^{-1}).
$$
It follows that if $L$ is a section of the bundle $Hom(TM,TM)$,
$$
\varphi_{\ast}\big((\nabla^{\ast}_{X}L)(Y)\big)=\big(\widetilde\nabla_{\varphi_{\ast}X}(\varphi_{\ast}\circ
L\circ \varphi_{\ast}^{-1})\big)(\varphi_{\ast}\circ
Y\circ\varphi^{-1}),
$$
where $\varphi_{\ast}\circ L\circ \varphi_{\ast}^{-1}(Z_{\widetilde
q})=\varphi_{\ast\,\varphi^{-1}(\widetilde q)}\circ
L_{\varphi^{-1}(\widetilde q)} \circ\varphi_{\ast\,\widetilde
q}^{-1}(Z_{\widetilde q})$ for every $\widetilde q\in\widetilde M$
and $Z_{\widetilde q}\in T_{\widetilde q}\widetilde M$. Hence
\begin{equation}\label{eq 2-1}
\big(\widetilde\nabla_{\varphi_{\ast}X}(\Phi\circ
S\circ\varphi^{-1})\big)(\varphi_{\ast}\circ
Y\circ\varphi^{-1})=\varphi_{\ast}\big((\nabla^{\ast}_{X}S)(Y)\big)
\end{equation}
for any section $S$ of ${\mathcal Z}$.  Recall that the Levi-Civita
connections $\nabla^{\ast}$ of the metric $\varphi^{\ast}\widetilde
g=e^{2f}g$ and $\nabla$ of $g$ are related by
\begin{equation}\label{relation}
\nabla^{\ast}_{X}Y=\nabla_{X}Y +X(f)Y + Y(f)X-g(X,Y)\nabla f
\end{equation}
where $\nabla f$ is the gradient of $f$ defined by $g(\nabla
f,Z)=Z(f)$.  Now, let $I\in{\mathcal Z}$, and let $S$ be a section
of ${\mathcal Z}$ in a neighnourhood of the point $p=\pi(I)$ such
that $S(p)=I$ and $\nabla S|_p=0$. Then, by (\ref{relation}), for
every $X\in T_pM$,
\begin{equation}\label{eq 2-2}
(\nabla^{\ast}_{X}S)(Y)=(IY)(f)-g(X,IY)\nabla f-Y(f)IX+g(X,Y)I\nabla
f.
\end{equation}
It is convenient to define an endomorphism $V_{I,X}$ of $T_pM$ by
$$
V_{I,X}(Y)=(IY)(f)X-g(X,IY)\nabla f-Y(f)IX+g(X,Y)I\nabla f.
$$
Then, it follows from (\ref{eq 2-1}) and (\ref{eq 2-2}) that the
second identity of (\ref{Jk-hol}) is equivalent to
\begin{equation}\label{eq 2}
(-1)^{k+1}I\circ V_{I,X}=V_{I,IX}.
\end{equation}
Clearly, the latter identity is satisfied if $k=1$.  For $k=2$, it
is equivalent to $V_{I,X}=0$, i.e.,
$$
(IY)(f)IX-g(X,Y)\nabla f+Y(f)X-g(X,IY)I\nabla f=0
$$
for every $I\in{\mathcal Z}$ and every $X,Y\in T_{\pi(I)}M$. Let
$p\in M$, and let $I\in {\mathcal Z}$ be a  complex structure on
$T_pM$. For any vector $Y\in T_pM$, take a non-zero vector $X\in
T_{p}M$ orthogonal to $Y$ and $IY$. Then the identity above implies
$(IY)(f)IX+Y(f)X=0$, hence $Y(f)=0$. It follows $f=const$.

\end{proof}

\smallskip

\noindent {\bf Remark}. The identity  $\Phi(I)=\varphi_{\ast}\circ
I\circ\varphi_{\ast}^{-1}$ follows from conformality of $\varphi$ if
$\varphi=Id$ and $\Phi$ is the map $\Psi$ defined above. This is not
the case in general.

\medskip

\noindent {\bf Notation}. Let $g_s$ and $\widetilde g_t$ be the
$1$-parameter families of Riemannian metrics on ${\mathcal Z}$ and
$\widetilde{\mathcal Z}$ corresponding to the metrics $g$ and
$\widetilde g$ on $M$. For fixed $s$ and $t$, let $D=D^{s}$ and
$\widetilde D=\widetilde D^{t}$ be the Levi-Civita connections of
the metrics  $g_s$ and $\widetilde g_t$.
\hfill $\Box$

\medskip

We are going to explore the problem of when the map $\Psi:({\mathcal
Z},g_s)\to (\widetilde{\mathcal Z},\widetilde g_t)$ is harmonic.

\smallskip

Let $\Psi^{\ast}T\widetilde{\mathcal Z}\to {\mathcal Z}$ be the
pull-back bundle of the bundle $T\widetilde{\mathcal Z}\to
\widetilde{\mathcal Z}$ under the map $\Psi:{\mathcal
Z}\to\widetilde{\mathcal Z}$. Then the differential
$\Psi_{\ast}:T{\mathcal Z}\to T\widetilde{\mathcal Z}$ can be
considered as a section of the vector bundle $Hom(T{\mathcal
Z},\Psi^{\ast}T\widetilde{\mathcal Z})$. Denote by ${\widetilde
D}^{\ast}$ the connection on the bundle
$\Psi^{\ast}T\widetilde{\mathcal Z}$ induced by the connection
$\widetilde D$ on $T\widetilde{\mathcal Z}$. The connections $D$ on
$T{\mathcal Z}$ and ${\widetilde D}^{\ast}$ on
$\Psi^{\ast}T\widetilde{\mathcal Z}$ induce a connection $\widehat
D$ on the bundle $Hom(T{\mathcal Z},\Psi^{\ast}T\widetilde{\mathcal
Z})$.

Recall that the second fundamental form of the map $\Psi$ is
$II_{\Psi}(X,Y)=(\widehat D_{X}\Psi_{\ast})(Y)$. It is a symmetric
$2$-form. Recall also that the map $\Psi:({\mathcal Z},g_s)\to
(\widetilde{\mathcal Z},\widetilde g_t)$ is harmonic if and only if
$$
Trace_{g_s}\widehat D\Psi_{\ast}=0.
$$
In order to compute the second fundamental form of the map $\Psi$ we
need the following lemma.

\begin{lemma}{\rm (\cite{D17, DM18})}\label{v-frame}
Let $K$ be a section of the twistor space $\widetilde{\mathcal Z}$
near a point $q\in M$. Then there exists a $\widetilde
g_t$-orthonormal frame of vertical vector fields
$\{V_{\alpha}:~\alpha=1,...,m^2-m\}$ on $\widetilde{\mathcal Z}$
defined in a neighbourhood of the point $K(q)$ such that

\noindent $(1)$ $\quad (\widetilde
D_{V_{\alpha}}V_{\beta})_{K(q)}=0$,~~ $\alpha,\beta=1,...,m^2-m$.

\noindent $(2)$ $\quad$ If $X$ is a vector field near the point $q$,
$[X^h,V_{\alpha}]_{K(q)}=0$.

\noindent $(3)$ $\quad$  $\widetilde\nabla_{X_q}(V_{\alpha}\circ
K)\perp {\mathcal V}_{K(q)}$

\end{lemma}

\begin{lemma}\label{sec fund form for hor}
Define a symmetric $2$-form $\Sigma$ on $M$ with values in $TM$ by
$$
\widetilde g(\Sigma(X,Y),Z)=\widetilde g(C^{-1}\circ
(\nabla_{Z}C)(X),Y),\quad X,Y,Z\in TM.
$$
Let $I\in{\mathcal Z}$, and  let $S$ be a section of ${\mathcal Z}$
such that $S(p)=I$, $p=\pi(I)$, and $\nabla S|_p=0$. Set
$J=\Psi(I)$.Then
$$
\begin{array}{c}
II_{\Psi}(X^h_I,Y^h_I)= -t\widetilde
R(\big(J\circ\widetilde\nabla_{X_p}(\Psi\circ
S)\big)^{\wedge})Y)^{\tilde h}_J-t\widetilde
R(\big(J\circ\widetilde\nabla_{Y_p}(\Psi\circ
S)\big)^{\wedge})X)^{\tilde h}_J\\[6pt]

+\displaystyle{\frac{1}{2}}\big(C^{-1}\circ(\nabla_{X}C)(Y)+C^{-1}\circ(\nabla_{Y}C)(X)-\Sigma(X,Y)\big)^{\tilde h}_J,\\[6pt]

+\displaystyle{\frac{1}{2}}\widetilde{\mathcal
V}_{J}\big(\widetilde\nabla^{2}_{XY}(\Psi\circ
S)+\widetilde\nabla^{2}_{YX}(\Psi\circ S)-\Psi\circ(\nabla^2_{XY}S)-\Psi\circ(\nabla^2_{YX}S)\big)\\
\end{array}
$$
where $\widetilde{\mathcal V}_J$ means ''the vertical component in
$T_J\widetilde{\mathcal Z}$'' and $\widetilde\nabla^2$, $\nabla^2$
stand for the second covariant derivatives,
$\nabla^2_{XY}S=\nabla_{X}\nabla_{Y}S-\nabla_{\nabla_{X}Y}S$, and
similarly for $\widetilde\nabla^{2}_{XY}(\Psi\circ S)$.
\end{lemma}

\begin{proof}
 Extend the
vector $Y\in T_pM$ to a vector field in a neighbourhood of $p$ such
that $\nabla Y|_p=0$.

Set  $J=\Psi(I)$ and $q=\widetilde \pi(J)$. Let
$\{V_{\alpha}:~\alpha=1,...,m^2-m\}$ be a frame of vertical vector
fields on $\widetilde{\mathcal Z}$ in a neighbourhood of the point
$J$ with the properties $(1)$ - $(3)$ stated in Lemma~\ref{v-frame}
with $K=\Psi\circ S$. Take a $\widetilde g$-orthonormal frame $
\tilde E_1,...,\tilde E_n$ in a neighbouhood of the point $q$ such
that $\widetilde{\nabla} \tilde E_i|_q=0$.

Then
$$
\Psi_{\ast}\circ Y^h=\sum\limits_{i=1}^n \widetilde
g_t(\Psi_{\ast}\circ Y^h,{\tilde E}^{\tilde h}_i\circ\Psi)({\tilde
E}^{\tilde h}
_i\circ\Psi)+\sum\limits_{\alpha=1}^{m^2-m}h_t(\Psi_{\ast}\circ
Y^h,V_{\alpha}\circ\Psi)(V_{\alpha}\circ\Psi).
$$
Hence, since $X^h_I=S_{\ast}X$,
$$
\begin{array}{r}
D^{\ast}_{X^h_I}(\Psi_{\ast}\circ Y^h)=D_{X}(\Psi_{\ast}\circ
Y^h\circ S)= (\sum\limits_{i=1}^n X_p\big(\widetilde
g_t(\Psi_{\ast}\circ Y^h\circ S,{\tilde E}^{\tilde
h}_i\circ\Psi\circ S)\big)({\tilde E}^{\tilde h}_i)_J\\[6pt]
+\sum\limits_{i=1}^n \widetilde g_t(\Psi_{\ast}(Y^h_I),({\tilde
E}^{\tilde
h}_i)_J)\widetilde D_{\Psi_{\ast}X^h_I}{\tilde E}_i^{\tilde h}\\[6pt]
+\sum\limits_{\alpha=1}^{m^2-m}X_p \big(\widetilde
g_t(\Psi_{\ast}\circ
Y^h\circ S,V_{\alpha}\circ\Psi\circ S)\big)(V_{\alpha})_J\\[6pt]
+\sum\limits_{\alpha=1}^{m^2-m}\widetilde g_t(\Psi_{\ast}(
Y^h_I),(V_{\alpha})_J)\widetilde D_{\Psi_{\ast}X^h_I}V_{\alpha}.
\end{array}
$$
Now, note that
\begin{equation}\label{no name}
\Psi_{\ast}\circ Y^h\circ
S=\Psi_{\ast}(S_{\ast}Y-\nabla_{Y}S)=Y^{\tilde h}\circ\Psi\circ
S+\widetilde\nabla_{Y}(\Psi\circ S)-\Psi_{\ast}(\nabla_{Y}S),
\end{equation}
where the last two terms are vertical tangent vectors of
$\widetilde{\mathcal Z}$. Then
$$
\sum\limits_{i=1}^n X_p\big(\widetilde g_t(\Psi_{\ast}\circ Y^h\circ
S,\tilde E^{\tilde h}_i\circ\Psi\circ S)\big)(\tilde E^{\tilde
h}_i)_J =\sum\limits_{i=1}^n X_p\big(\widetilde g(Y,\tilde
E_i)\big)(\tilde E^{\tilde h}_i)_J=(\widetilde\nabla_{X}Y)^{\tilde
h}_J.
$$
Also, by Lemma~\ref{Levi-Civita},
$$
\begin{array}{c}
\sum\limits_{i=1}^n \widetilde g_t(\Psi_{\ast}(Y^h_I),(\tilde
E^{\tilde
h}_i)_J)\widetilde D_{\Psi_{\ast}X^h_I}\tilde E_i^{\tilde h}\\[6pt]
=\sum\limits_{i=1}^n \widetilde g_t(Y^{\tilde h}_J,(\tilde E^{\tilde
h}_i)_J)\big[\widetilde D_{X^{\tilde h}_J}\tilde E^{\tilde
h}_i+\widetilde D_{\widetilde\nabla_{X_p}(\Psi\circ S)}\tilde E^{\tilde h}_i\big]\\[6pt]

=\displaystyle{\sum\limits_{i=1}^n \widetilde g(Y_p,(\tilde
E_i)_p)\big[\frac{1}{2}\widetilde R(X,\tilde E_i)J-t(\widetilde
R(\big(J\circ\widetilde\nabla_{X_p}(\Psi\circ
S)\big)^{\wedge})\tilde E_i)^{\tilde h}_J\big]} \\[10pt]

=\displaystyle{\frac{1}{2}\widetilde R(X,Y)J-t(\widetilde
R(\big(J\circ\widetilde \nabla_{X_p}(\Psi\circ
S)\big)^{\wedge})Y)^{\tilde h}_J},
\end{array}
$$
where  in the first summand  $\widetilde R$ stands for the curvature
of the connection $\widetilde\nabla$ on the bundle $A(TM)$, while in
the second summand it stands for the curvature of $TM$.

By (\ref{G-H}), (\ref{diff on ver}), and (\ref{no name}),

$$
\begin{array}{c}
\sum\limits_{\alpha=1}^{m^2-m}X_p \big(\widetilde
g_t(\Psi_{\ast}\circ Y^h\circ
S,V_{\alpha}\circ\Psi\circ S)\big)(V_{\alpha})_J\\[6pt]

=\sum\limits_{\alpha=1}^{m^2-m}tX_p \big(\widetilde
G(\widetilde\nabla_{Y}(\Psi\circ
S),V_{\alpha}\circ\Psi\circ S)\big)(V_{\alpha})_J \\[6pt]

-\sum\limits_{\alpha=1}^{m^2-m}tX_p \big(G(\nabla_{Y}S,Q\circ
(V_{\alpha}\circ\Psi\circ S)\circ Q^{-1}\big)(V_{\alpha})_J\\[6pt]

=\sum\limits_{\alpha=1}^{m^2-m} \widetilde
g_t(\widetilde\nabla_{X}\widetilde\nabla_{Y}(\Psi\circ
S),V_{\alpha}\circ\Psi\circ S)(V_{\alpha})_J\\[6pt]

-\sum\limits_{\alpha=1}^{m^2-m}tG(\nabla_{X}\nabla_{Y}S,Q\circ
(V_{\alpha}\circ\Psi\circ S)\circ Q^{-1}\big)(V_{\alpha})_J

\end{array}
$$
since $\widetilde\nabla_{X_p}(V_{\alpha}\circ\Psi\circ S)\perp
{\mathcal V}_{J}$ and $\nabla_{X_p}S=0$.

Therefore
$$
\sum\limits_{\alpha=1}^{m^2-m}X_p \big(\widetilde
g_t(\Psi_{\ast}\circ Y^h\circ S,V_{\alpha}\circ\Psi\circ
S)\big)(V_{\alpha})_J=\widetilde{\mathcal
V}_{J}\big(\widetilde\nabla_{X}\widetilde\nabla_{Y}(\Psi\circ
S)-Q^{-1}\circ (\nabla_{X}\nabla_{Y}S)\circ Q\big)
$$
where $\widetilde{\mathcal V}$ on the right-hand side stands for the
vertical component in $T\widetilde{\mathcal Z}$.

By Lemmas~\ref{v-frame} and \ref{Levi-Civita},
$$
\widetilde D_{\Psi_{\ast}X^h_I}V_{\alpha}=\widetilde D_{X^{\tilde
h}_J}V_{\alpha}=\widetilde D_{(V_{\alpha})_J}X^{\tilde
h}=-t\big(\widetilde R((J\circ
(V_{\alpha})_J)^{\wedge})X\big)^{\tilde h}_J.
$$
Hence
$$
\sum\limits_{\alpha=1}^{m^2-m}\widetilde g_t(\Psi_{\ast}\circ
Y^h_I,(V_{\alpha})_J)\widetilde
D_{\Psi_{\ast}X^h_I}V_{\alpha}=-t\widetilde
R(\big(J\circ\widetilde\nabla_{Y_p}(\Psi\circ
S)\big)^{\wedge})X)^{\tilde h}_J).
$$
Finally, note that by Lemma~\ref{Levi-Civita},
$$
\Psi_{\ast}(D_{X^h_I}Y^h)=\Psi_{\ast}(\frac{1}{2}R_p(X,Y)I)
=\frac{1}{2}Q^{-1}\circ R_p(X,Y)I\circ Q.
$$
It follows that
$$
\begin{array}{c}
II_{\Psi}(X^h_I,Y^h_I)= -t\widetilde
R(\big(J\circ\widetilde\nabla_{X_p}(\Psi\circ
S)\big)^{\wedge})Y)^{\tilde h}_J-t\widetilde
R(\big(J\circ\widetilde\nabla_{Y_p}(\Psi\circ
S)\big)^{\wedge})X)^{\tilde h}_J\\[6pt]

+\displaystyle{\frac{1}{2}}\widetilde{\mathcal
V}_{J}\big(\widetilde\nabla^{2}_{XY}(\Psi\circ
S)+\widetilde\nabla^{2}_{YX}(\Psi\circ S)-Q^{-1}\circ
\nabla^2_{XY}S\circ
Q-Q^{-1}\circ \nabla^2_{YX}S\circ Q\big)\\[6pt]

+(\widetilde\nabla_{X}Y)^{\tilde h}_J.
\end{array}
$$

The Koszul formula for the Levi-Civita connection and identity
(\ref{C}) imply the identity
$$
\begin{array}{c}
2\widetilde g(\widetilde\nabla_{X}Y,Z)=2\widetilde g(\nabla_{X}Y,Z)\\[6pt]
+\widetilde g(C^{-1}(\nabla_{X}C)(Y),Z)+\widetilde
g(C^{-1}(\nabla_{Y}C)(X),Z)-\widetilde g(C^{-1}(\nabla_{Z}C)(X),Y)
\end{array}
$$
for any $X,Y,Z\in TM$. Now, the lemma follows from the last two
identities.
\smallskip
\end{proof}

\begin{lemma}\label{sec fund form for ver}
If $U'$ and $U''$ are vertical vectors at a point $I\in {\mathcal
Z}$,
$$
II_{\Psi}(U',U'')=0.
$$
\end{lemma}

\begin{proof}
Let $\{V_{\alpha}:~\alpha=1,...,m^2-m\}$ be a $\widetilde
g_t$-orthonormal frame of vertical vector fields on
$\widetilde{\mathcal Z}$ in a neighbourhood of the point $J=\Psi(I)$
such that $\widetilde D_{V_{\alpha}(J)}V_{\beta}=0$. The map $\Psi$
sends the fibre ${\mathcal Z}_p$, $p=\pi(I)$, of ${\mathcal Z}$
through $I$ onto the fibre $\widetilde{\mathcal Z}_q$,
$q=\widetilde\pi(J)$, of $\widetilde{\mathcal Z}$ through $J$. By
(\ref{G-H}), $\Psi|{\mathcal Z}_p\to \widetilde{\mathcal Z}_q$ is an
isomorphism with respect to the restrictions of the metrics
$\frac{1}{s}g_s$ and $\frac{1}{t}\widetilde g_t$ to ${\mathcal Z}_p$
and $\widetilde{\mathcal Z}_q$, respectively. Hence, $\Psi|{\mathcal
Z}_p$ preserves the Levi-Civita connections of these metrics on
${\mathcal Z}_p$ and $ \widetilde{\mathcal Z}_q$. These connections
coincide with the Levi-Civita connections of the metrics
$g_s|{\mathcal Z}_p$ and $\widetilde g_t|\widetilde{\mathcal Z}_q$.
Note also that the fibres of ${\mathcal Z}$ and $\widetilde{\mathcal
Z}$ are totally geodesic submanifolds. Hence $\Psi|{\mathcal Z}_p$
preserves the restrictions of the connections $D$ to ${\mathcal
Z}_p$ and $\widetilde D$ to $\widetilde{\mathcal Z}_q$. Therefore $
U_{\alpha}=\frac{t}{\sqrt s}\Psi_{\ast}^{-1}(V_{\alpha}\circ \Psi)$,
$\alpha=1,...,m^2-m, $ constitute a $g_s$-orthonormal frame of
vertical vector fields on ${\mathcal Z}$ in a neighbourhood of $I$
such that $D_{U_{\alpha}(I)}U_{\beta}=0$. Then, at the point $I$,
$$
II_{\Psi}(U_{\alpha},U_{\beta})={D^{\ast}}_{U_{\alpha}(I)}(V_{\beta}\circ\Psi)-\Psi_{\ast}(D_{U_{\alpha}(I)}U_{\beta})=\widetilde
D_{V_{\alpha}(J)}V_{\beta}=0.
$$
This proves the lemma.
\end{proof}

\smallskip

\noindent {\it Assumption}. The most interesting case for us is that
of conformally related metrics, i.e.,  $\widetilde g=e^{2f}g$ for a
smooth function $f$.  Henceforth we assume that $C=e^{2f}Id$, hence
$\Psi$ is the identity map.

\smallskip

Clearly, if  $\widetilde g=e^{2f}g$,  ${\mathcal
Z}=\widetilde{\mathcal Z}$ as smooth bundles over $M$. Also, they
have the same Atiyah-Hitchin-Singer almost complex structure. But,
endowed with the metrics $g_s$ and $\widetilde g_t$,  ${\mathcal Z}$
and $\widetilde{\mathcal Z}$  are different as Riemannian manifolds.

In the case $C=e^{2f}Id$, Lemmas~\ref{sec fund form for hor} and
\ref{sec fund form for ver} take the following form.

\begin{lemma}\label{case Id}
Under the notation in  Lemmas~\ref{sec fund form for hor} and
\ref{sec fund form for ver},
$$
{\mathcal V}
II_{Id}(X^h_I,Y^h_I)=\displaystyle{\frac{1}{2}}{\mathcal
V}\big(\widetilde\nabla^{2}_{XY} S+\widetilde\nabla^{2}_{YX}S-
\nabla^2_{XY}S- \nabla^2_{YX}S\big),
$$
$$
\begin{array}{c}
{\mathcal H} II_{Id}(X^h_I,Y^h_I)=-t\widetilde
R(\big(I\circ\widetilde\nabla_{X_p}S\big)^{\wedge})Y)^{\tilde
h}_I-t\widetilde
R(\big(I\circ\widetilde\nabla_{Y_p}S\big)^{\wedge})X)^{\tilde
h}_I\\[6pt]

\hspace{2cm}+(X(f)Y+Y(f)X-g(X,Y)\nabla f)^{\tilde h}_I.
\end{array}
$$

\end{lemma}

\begin{prop}\label {Id harm} Suppose that $g$ and $\widetilde g$ are
conformal metrics on a manifold $M$, $\widetilde g=e^{2f}g$. The
identity map $({\mathcal Z},g_s)\to (\widetilde{\mathcal
Z},\widetilde g_t)$ is harmonic if and only if the metrics $g$ and
$\widetilde g$ are homothetic, $f=const$.

\end{prop}

\begin{proof}

Recall that in this case
\begin{equation}\label{LC}
\widetilde\nabla_{X}Y=\nabla_{X}Y +X(f)Y + Y(f)X-g(X,Y)\nabla f
\end{equation}
Hence, for any endomorphism $L$ of $TM$,
\begin{equation}\label{L}
(\widetilde\nabla_{X}
L)(Y)=(\nabla_{X}L)(Y)+(LY)(f)X-Y(f)L(X)-g(X,L(Y))\nabla
f+g(X,Y)L(\nabla f).
\end{equation}
In the next  computation, tt is convenient to use the isomorphism
$A(TM)\cong \Lambda^2TM$ defined via identity (\ref{wedge-isom}). If
$\sigma\in\Lambda^2TM$, the corresponding skew-symmetric
endomorphism will be denoted by $\sigma^{\vee}$; it is determined by
\begin{equation}\label{endo}
g(\sigma^{\vee}(X),Y)=2g(\sigma,X\wedge Y),\quad X,Y\in TM.
\end{equation}

With this notation, if $L$ is a skew-symmetric endomorphism of $TM$,
identity (\ref{L}) can be written as
\begin{equation}\label{L-skew}
\widetilde\nabla_{X} L=\nabla_{X}L-(L(\nabla f)\wedge X+\nabla
f\wedge L(X))^{\vee}.
\end{equation}

Now, let $A,B$ be vector fields on $M$  and $X\in TM$. Then
\begin{equation}\label{vee}
\begin{array}{c}
(A\wedge B)^{\vee}(X)=g(A,X)B-g(B,X)A,\\[6pt]
 \nabla_{X}(A\wedge
B)^{\vee}=(\nabla_{X}A\wedge B+A\wedge \nabla_{X}B)^{\vee}.
\end{array}
\end{equation}
The first of the these identities is a direct consequence from
(\ref{endo}), the proof of the second one is standard. Computing
$\widetilde\nabla_{X}A\wedge B+A\wedge\widetilde\nabla_{X}B$ by
means of (\ref{LC}), and applying (\ref{vee}) and (\ref{L-skew}) to
compute $\widetilde\nabla_{X}(A\wedge B)^{\vee}$, we see also that
\begin{equation}\label{tilde-vee}
\widetilde\nabla (A\wedge B)^{\vee}=\big(\widetilde\nabla_{X}\wedge
B+A\wedge\widetilde\nabla_{X}B\big)^{\vee}-2X(f)(A\wedge B)^{\vee}.
\end{equation}

Now, we discuss  the vertical part of the second fundamental form of
the map $Id: ({\mathcal Z},g_s)\to ({\mathcal Z},\widetilde g_t)$.

Let $X,Y\in T_pM$. Extend $Y$ to a vector field in a neighbourhood
of $p$ such that $\nabla Y|_p=0$. Also, let $S$ be a section of
${\mathcal Z}$ as in Lemma~\ref{sec fund form for hor}, so $\nabla
S|_p=0$. Then an easy computation making use of (\ref{tilde-vee})
gives
\begin{equation}\label{rel sec der}
\begin{array}{l}
\widetilde\nabla^{2}_{XY}S=\widetilde\nabla_{X}\widetilde\nabla_{Y}S-\widetilde\nabla_{\widetilde\nabla_{X}Y}S\\[8pt]

=\widetilde\nabla_{X}\nabla_{Y}S-\big(\widetilde\nabla_{X}S(\nabla
f)\wedge Y+\widetilde\nabla_{X}(\nabla f)\wedge S(Y)+\nabla f\wedge
\widetilde\nabla_{X}S(Y)\big)^{\vee}\\[6pt]

-2X(f)(S(\nabla f)\wedge Y +(\nabla f)\wedge S(Y)-\nabla f\wedge S(\widetilde\nabla_{X}Y)\big)^{\vee}\\[8pt]

=\nabla^2_{XY}S +2X(f)(S(\nabla f)\wedge Y +(\nabla f)\wedge S(Y))^{\vee}\\[6pt]
-\big(S(\nabla_{X}(\nabla f))+X(f)S(\nabla f)+S(\nabla
f)(f)X+S(X)(f)\nabla f\big)\wedge Y\big)^{\vee}\\[6pt]

-\big(\nabla_{X}(\nabla f)+X(f)\nabla f+(\nabla f)(f)X-X(f)\nabla f\big)\wedge S(Y)\big)^{\vee}\\[6pt]

-\big((\nabla f)\wedge S(Y)(f)X-Y(f)S(X)+g(X,Y)S(\nabla
f)\big)^{\vee}.

\end{array}
\end{equation}

Next, note that if $A,B\in TM$, for every $X,Y\in TM$,
$$
\begin{array}{c}
g((S\circ (A\wedge B)^{\vee}\circ S)^{\wedge},X\wedge
Y)=-\displaystyle{\frac{1}{2}}g((A\wedge B)^{\vee}(S(X)),S(Y))\\[6pt]
=-g(A\wedge B,S(X)\wedge S(Y))=-g(S(A)\wedge S(B),X\wedge Y).
\end{array}
$$
Hence
\begin{equation}\label{S-S}
S\circ (A\wedge B)^{\vee}\circ S=-\big(S(A)\wedge S(B)\big)^{\vee}.
\end{equation}
The latter identity and (\ref{rel sec der}) imply
\begin{equation}\label{v-proj}
\begin{array}{l}
\widetilde\nabla^{2}_{XY}S+S\circ \widetilde\nabla^{2}_{XY}S\circ
S=\nabla^{2}_{XY}S+S\circ \nabla^{2}_{XY}S\circ S \\[8pt]

-S(\nabla f)(f)\big(X\wedge Y-S(X)\wedge S(Y)\big)^{\vee}-S(X)(f)\big(\nabla f\wedge Y-S(\nabla f)\wedge S(Y)\big)^{\vee}\\[6pt]
-(\nabla f)(f)\big(X\wedge S(Y)+S(X)\wedge
Y\big)^{\vee}+X(f)\big(\nabla f\wedge S(Y)+S(\nabla f)\wedge
Y\big)^{\vee}.
\end{array}
\end{equation}
By (\ref{S-S}) and (\ref{proj}), every endomorphism of the form
$\big(A\wedge B-S(A)\wedge S(B)\big)^{\vee}$ is, in fact, a vertical
tangent vector of ${\mathcal Z}$. This remark and (\ref{v-proj}),
imply
$$
\begin{array}{c}
{\mathcal V}\widetilde\nabla^{2}_{XX}S -{\mathcal
V}\nabla^{2}_{XX}S\\[6pt]

=-g(S(X),\nabla f)\big(\nabla f\wedge X-S(\nabla f)\wedge
S(X)\big)^{\vee}+g(X,\nabla f)\big(\nabla f\wedge S(X)+S(\nabla
f)\wedge X\big)^{\vee}.
\end{array}
$$
Let $E_1,...,E_n$ be an orthonormal basis of the tangent space
$T_pM$. Then, for any $A,B\in T_pM$,
$$
\sum\limits_{i=1}^n g(S(E_i),\nabla f)g(\nabla f\wedge E_i-S(\nabla
f)\wedge S(E_i),A\wedge B)=0.
$$
Hence
$$
Trace\{X\to g(S(X),\nabla f)\big(\nabla f\wedge X-S(\nabla f)\wedge
S(X)\big)^{\vee}\}=0.
$$
Also,
$$
Trace\{X\to g(X,\nabla f)\big(\nabla f\wedge S(X)+S(\nabla f)\wedge
X\big)^{\vee}\}=0.
$$
Therefore
$$
Trace\,{\mathcal V}(\widetilde\nabla^{2}_{XX}S-\nabla^{2}_{XX}S)=0.
$$
Next, we turn to the horizontal part of the second fundamental form
$II_{Id}$. Note that, for $A,B,X,Y\in TM$,
$$
\begin{array}{c}
g(\big(S\circ (S(A)\wedge B+A\wedge S(B))^{\vee}\big)^{\wedge},X\wedge Y)\\[6pt]

=-\displaystyle{\frac{1}{2}}g((S(A)\wedge B+A\wedge S(B))^{\vee}(X),S(Y))\\[6pt]

=-g(S(A)\wedge B+A\wedge S(B),X\wedge S(Y))\\[6pt]

=g(S(A)\wedge S(B)-A\wedge B,X\wedge Y).

\end{array}
$$
Hence
$$
\big(S\circ (S(A)\wedge B+A\wedge
S(B))^{\vee}\big)^{\wedge}=S(A)\wedge S(B)-A\wedge B.
$$
This identity and (\ref {L-skew}) imply
$$
\big(I\circ\widetilde\nabla_{X_p}S\big)^{\wedge}=\nabla f\wedge
X-I(\nabla f)\wedge I(X),\quad I\in{\mathcal Z},\quad p=\pi(I).
$$
Now, it follows from Lemma~\ref{case Id}, that the identity map
$({\mathcal Z}, g_s)\to ({\mathcal Z},\widetilde g_t)$ is harmonic
if and only if for every $I\in {\mathcal Z}$ and every $X\in
T_{\pi(I)}M$,
\begin{equation}\label{Id harmonic}
\begin{array}{c}
2te^{2f}\widetilde\rho(\nabla f,X)+2tTrace_{g}\{T_{\pi(I)}M\ni Z\to
\widetilde g(\widetilde R(I(\nabla f)\wedge I(Z))Z,X)\}\\[6pt]
-(n-2)\widetilde g(\nabla f,X)=0,
\end{array}
\end{equation}
where $\widetilde\rho$ is the Ricci tensor of the metric $\widetilde
g$.

Let $E_1,...,E_n$ be a $g$-orthonormal basis of $T_{\pi(I)}M$ such
that $E_{2k}=I(E_{2k-1})$, $k=1,...,m$.  By the algebraic Bianchi
identity
$$
\begin{array}{c}
Trace_{g}\{T_{\pi(I)}M\ni Z\to \widetilde g(\widetilde R(I(\nabla
f)\wedge I(Z))Z,X) \\[6pt]
=\displaystyle{\frac{1}{2}}\sum\limits_{i=1}^n\widetilde
g(\widetilde R(E_i\wedge I(E_i))I(\nabla f),X).
\end{array}
$$

Suppose that the map $Id: ({\mathcal Z}, g_s)\to ({\mathcal
Z},\widetilde g_t)$ is harmonic.

Fix a point $p\in M$. Identity (\ref{Id harmonic}) is trivially
satisfied at the points where $\nabla f$ vanishes, so suppose
$(\nabla f)(p)\neq 0$  at a point $p$ of $M$. Set $E_1=\frac{(\nabla
f)(p)}{||(\nabla f)(p)||_{g}}$. Let
$e=(E_1,E_2,E_3,E_4,...,E_{n-1},E_{n})$ be any oriented
$g$-orthonormal basis of $T_pM$. Define a complex structure $I$ on
$T_pM$ by $I(E_{2k-1})=E_{2k}$, $I(E_{2k})=-E_{2k-1}$, $k=1,...,m$,
$2m=n$. Clearly $I\in{\mathcal Z}$, hence
\begin{equation}\label{main eq}
\begin{array}{c}
2t\widetilde g(\widetilde R(E_1\wedge E_2+E_3\wedge
E_4+\cdots +E_{2m-1}\wedge E_{2m})E_2,X)\\[6pt]
+2te^{2f}\widetilde\rho(E_1,X) -(n-2)\widetilde
g(E_1,X)=0~~{\rm{for~every}}~~X\in T_{\pi(I)}M.
\end{array}
\end{equation}
The $g$-orthonormal basis $e'=(E_1,-E_2,-E_3,E_4,E_5,...,E_n)$  is
positively oriented. Applying identity (\ref{main eq}) for this
basis and comparing with the same  identity for the bases $e$, we
see that for every $X\in T_pM$,
\begin{equation}\label{eq 01}
\begin{array}{c}
2t\widetilde g(\widetilde R(E_1\wedge E_2)E_2,X)
+2t\widetilde g(\widetilde R(E_3\wedge E_4)E_2,X)\\[6pt]
+2te^{2f}\widetilde\rho(E_1,X)-(n-2)\widetilde g(E_1,X)=0,
\end{array}
\end{equation}
and if $n\geq 6$,
\begin{equation}\label{eq 02}
\widetilde g(\widetilde R(E_5\wedge E_6+\cdots +E_{n-1}\wedge
E_{n})E_2,X)=0.
\end{equation}

We discuss first the case $n=4$. Add the three identities obtained
by writing (\ref{eq 01}) for the bases $(E_1,E_2,E_3,E_4)$,
$(E_1,E_3,E_4,E_2)$ and $(E_1,E_4,E_2,E_3)$. Taking into account the
algebraic Bianchi identity, we get
\begin{equation}\label{case 1-eq 01}
2te^{2f}\widetilde \rho(E_1,X)-3\widetilde g(E_1,X)=0.
\end{equation}
It follows that for every oriented $g$-orthonormal basis
$(E_1,...,E_4)$ of $T_pM$,
\begin{equation}\label{case 1-eq 02}
2t\widetilde R(E_1\wedge E_2+E_3\wedge E_4)E_2=-E_1.
\end{equation}
Let $s_1=s_1^{+}, s_2=s_2^{+}, s_3=s_3^{+}$ be the basis of
$\Lambda^2T_pM$ defined by (\ref{s-basis}) for the basis
$e=(E_1,E_2,E_3,E_4)$ . Let $e'=(E_1',E_2',E_3',E_4')$ be the
oriented $g$-orthonormal basis of $T_pM$ obtained from $e$ by the
action of a matrix $A\in SO(4)$. Denote by $s_1',s_2',s_3'$ the
$2$-vectors corresponding to the basis $e'$ via (\ref{s-basis}). It
is well-known that every matrix $A$ in $SO(4)$ can be represented as
the product $A=A_1A_2$ of two $SO(4)$-matrices of the following
types
\begin{equation}\label{isoclinic}
A_1=\left(\begin{array}{cccc}
                   a & -b & -c & -d \\
                   b & a & -d & c \\
                   c & d & a & -b \\
                   d & -c & b & a \\
                  \end{array} \right),\quad
A_2=\left(\begin{array}{cccc}
                   p & -q & -r & -s \\
                   q & p & s & -r \\
                   r & -s & p & q \\
                   s & r & -q & p \\
                  \end{array} \right),
\end{equation}
where $a,...,d,p,...,s$ are real numbers with $a^2+b^2+c^2+d^2=1$,
$p^2+q^2+r^2+s^2=1$, the isoclinic representation found by L. van
Elfrikhof \cite{E}(1897). Note that the representation $A=A_1A_2$
corresponds to the well-known isomorphism $SO(4)\cong (Sp(1)\times
Sp(1))/{\mathbb Z}_2$. The transformation $A$ preserves the vector
$E_1$ exactly when $p=a,q=-b,r=-c,a=-d$. In this case $A$ has the
form
\begin{equation}\label{A}
A=\left(\begin{array}{cccc}
                   1 & 0 & 0 & 0 \\
                   0 &  &    &  \\
                   0 &  & \widehat A &  \\
                   0 &  &   & \\
                  \end{array} \right)
\end{equation}
where $\widehat A$ is the $SO(3)$-matrix
$$
\widehat A=\left( \begin{array}{cccc}
                                       a^2+b^2-(c^2+d^2) & -2(ad-bc) & 2(ac+bd)  \\
                                       2(ad+bc) & a^2+c^2-(b^2+d^2) & -2(ab-cd) \\
                                       -2(ac-bd) & 2(ab+cd) & a^2+d^2-(b^2+c^2)\\
                                                                            \end{array}
                                   \right).
$$
This matrix represents the standard action of the unit quaternion
$a+ib+jc+kd$ on ${\mathbb R}^3$. Writing identity (\ref{case 1-eq
02}) for the basis $(E_1',...,E_4')$, we get a polynomial identity
in $a,b,c,d$ which implies the identities
$$
\begin{array}{c}
\widetilde R(s_1)E_2=\widetilde R(s_2)E_3=\widetilde
R(s_3)E_4=-\displaystyle{\frac{1}{2t}}E_1, \\[8pt]
\widetilde R(s_1)E_3+\widetilde R(s_2)E_2=0,\quad \widetilde
R(s_1)E_4+\widetilde R(s_3)E_2=0,\quad \widetilde R(s_
2)E_4+\widetilde R(s_3)E_3=0.
\end{array}
$$
These identities imply
\begin{equation}\label{1-1}
\begin{array}{c}
\widetilde g(\widetilde{\mathcal R}(s_1),s_1^{\pm})=\widetilde
g(\widetilde R(s_1)E_1,E_2)\pm \widetilde
g(\widetilde R(s_1)E_3,E_4)\\[6pt]

=-\widetilde g(\widetilde R(s_1)E_2,E_1)\pm \widetilde g(\widetilde
R(s_2)E_4,E_2)\\[6pt]

=\widetilde g(\widetilde R(s_2)E_1,E_3)\pm \widetilde g(\widetilde
R(s_2)E_4,E_2)=\widetilde g(\widetilde{\mathcal R}(s_2),s_2^{\pm}).
\end{array}
\end{equation}
Also,
$$
\begin{array}{c}
\widetilde g(\widetilde{\mathcal R}(s_1),s_2^{\pm})=\widetilde
g(\widetilde R(s_2)E_2,E_1)\pm \widetilde g(\widetilde
R(s_1)E_4,E_2)\\[6pt] =-\widetilde g(\widetilde{\mathcal
R}(s_2),s_1)=-\widetilde g(\widetilde{\mathcal R}(s_2),s_1^{-})
\end{array}
$$
since $\widetilde g(\widetilde R(s_1)E_4,E_2)=\widetilde
g(\widetilde R(s_2)E_3,E_4)=0$. It follows
\begin{equation}\label{1-2}
\widetilde g(\widetilde{\mathcal R}(s_1),s_2^{\pm})=\widetilde
g(\widetilde{\mathcal R}(s_2),s_1^{-})=0.
\end{equation}
Similarly,
\begin{equation}\label{1-3}
\widetilde g(\widetilde{\mathcal R}(s_1),s_3^{\pm})=\widetilde
g(\widetilde{\mathcal R}(s_3),s_1^{-})=0.
\end{equation}
Identities (\ref{1-1}), (\ref{1-2}) and (\ref{1-3}) hold also for
$(s_1,s_2,s_3)$ replaced with $(s_2,s_3,s_1)$ and $(s_3,s_1,s_2)$.
It follows
$$
\widetilde g(\widetilde{\mathcal W}_{+}(s_i),s_j)=0\quad
{\rm{if}}~~i\neq j.
$$
Moreover,
$$
\widetilde g(\widetilde{\mathcal W}_{+}(s_1),s_1)=\widetilde
g(\widetilde{\mathcal W}_{+}(s_2),s_2)=\widetilde
g(\widetilde{\mathcal W}_{+}(s_3),s_3).
$$
Since $Trace\,\widetilde{\mathcal W}_{+}=0$, we get
$$
\widetilde g(\widetilde{\mathcal W}_{+}(s_i),s_i)=0.
$$
Thus, $\widetilde{\mathcal W}_{+}=0$.  We also have
$$
\begin{array}{c}
\widetilde g(\widetilde{\mathcal B}(s_i),s_j^{-})=0\quad
{\rm{if}}~~i\neq j,\\[6pt]

\widetilde g(\widetilde{\mathcal B}(s_1),s_1^{-})=\widetilde
g(\widetilde{\mathcal B}(s_2),s_2^{-})=\widetilde
g(\widetilde{\mathcal B}(s_3),s_3^{-}).
\end{array}
$$
Therefore  $\widetilde{\mathcal B}=0$.

As a conclusion, if identity (\ref{Id harmonic}) is satisfied,
setting $M'=\{p\in M: (\nabla f)(p)\neq 0\}$, the Riemannian
manifold $(M',\widetilde g)$ is Einstein and anti-self-dual. Then,
on $M'$,
$$
\frac{1}{2t}e^{2f}=\widetilde g(\widetilde
R(s_1)E_1,E_2)=\frac{1}{2}\widetilde g(\widetilde{\mathcal
R}(s_1),s_1+s_1^{-})=\frac{\tau}{12}\widetilde
g(s_1,s_1)=\frac{\tau}{12}e^{4f}.
$$
It follows that $f=const$ on $M'$, a contradiction. Thus $\nabla
f=0$ on $M$, i.e., $f=const$.

\smallskip

Now, suppose $n\geq 6$, and apply identity (\ref{eq 01}) for the
basis $e''=(E_1,E_2,-E_3,E_4,\\-E_5,E_6,E_7,E_8,...,E_{n-1},E_n)$.
The identity obtained, (\ref{eq 01}) and (\ref{eq 02}) imply
\begin{equation}\label{case 2-eq 1}
2t\widetilde g(\widetilde R(E_1\wedge E_2)E_2,X)\\[6pt]
+2te^{2f}\widetilde\rho(E_1,X)-(n-2)\widetilde g(E_1,X)=0
\end{equation}
and
\begin{equation}\label{case 2-eq 2}
\widetilde R(E_3\wedge E_4)E_2=0.
\end{equation}
Identity (\ref{case 2-eq 1}) holds for every unit vector $E_2$
orthogonal to $E_1$, in particular for $E_2$ replaced with
$E_3,...,E_n$. Adding the corresponding $(n-1)$ identities, we get
$$
2te^{2f}\widetilde \rho(E_1,X)-(n-1)\widetilde g(E_1,X)=0.
$$
The latter identity and  (\ref{case 2-eq 1}) imply
$$
2t\widetilde g(\widetilde R(E_1\wedge E_2)E_2,X)=-\widetilde
g(E_1,X).
$$
It follows that
\begin{equation}\label{case 2-eq 1-2}
2t\widetilde R(E_1,E_i)E_i=-E_1\quad {\rm{for}}~~i\geq 2.
\end{equation}
Also, by (\ref{case 2-eq 2}),
\begin{equation}\label{case 2-eq 2-2}
\widetilde R(E_i,E_j)E_k=0\quad {\rm{for}}~~i\neq j,\, j\neq k,\,
k\neq i,~ i,j,k\geq 2.
\end{equation}
As above, set $e=(E_1,...,E_n)$, and let $e'=(E_1',...,E_n')$ be the
oriented $g$-orthonormal basis of $T_pM$ such that $E_k'=E_k$ for
$k\geq 5$ and $(E_1',E_2',E_3',E_4')$ is obtained from
$(E_1,E_2,E_3,E_4)$ by the action of the matrix (\ref{A}). Identity
(\ref{case 2-eq 2}) applied for the basis $e'$ yields a polynomial
identity in $a,b,c,d$, and an easy but tedious computation gives
$$
\begin{array}{c}
\widetilde R(E_2,E_4)E_2=\widetilde R(E_3,E_4)E_3,\quad \widetilde
R(E_2,E_3)E_2=0,\\[6pt]

\widetilde R(E_2,E_3)E_3=\widetilde R(E_2,E_4)E_4,\quad \widetilde
R(E_3,E_4)E_4=0.
\end{array}
$$
It follows
\begin{equation}\label{case 2-eq 3}
\widetilde R(E_i,E_j)E_i=\widetilde R(E_i,E_j)E_j=0 \quad
{\rm{for}}~~i,j\geq 2.
\end{equation}
Writing identity (\ref{case 2-eq 1-2}) for the basis $e'$, we see
that
\begin{equation}\label{case 2-eq 4}
R(E_1,E_i)E_j+R(E_1,E_j)E_i=0 \quad {\rm{for}}~~  i,j\geq 2,~i\neq
j.
\end{equation}
Identity (\ref{case 2-eq 1-2}) implies
\begin{equation}\label{Ricci 1-1/i}
\widetilde\rho(E_1,E_1)=\frac{n-1}{2te^{2f}}\widetilde
g(E_1,E_1),\quad \widetilde\rho(E_1,E_i)=0\quad {\rm{for}}~~i\geq 2.
\end{equation}
It follows from (\ref{case 2-eq 1-2})  and (\ref{case 2-eq 3}) that
for $i\geq 2$
\begin{equation}\label{Ricci-i-i}
\begin{array}{c}
\widetilde\rho(E_i,E_i)=e^{-2f}\widetilde g(\widetilde
R(E_i,E_1)E_i,E_1)+e^{-2f}\sum\limits_{j=2}^n\widetilde g(\widetilde
R(E_i,E_j)E_i,E_j)\\[6pt]
=\displaystyle{\frac{1}{2te^{2f}}}\widetilde
g(E_1,E_1)=\displaystyle{\frac{1}{2te^{2f}}}\widetilde g(E_i,E_i).
\end{array}
\end{equation}
By (\ref{case 2-eq 2-2}), for $i,j\geq 2$, $i\neq j$,
$$
\widetilde\rho(E_i,E_j)=e^{-2f}\widetilde g(\widetilde
R(E_i,E_1)E_j,E_1).
$$
Hence, by (\ref{case 2-eq 4}),
$$
\begin{array}{c}
2\widetilde\rho(E_i,E_j)=\widetilde\rho(E_i,E_j)+\widetilde\rho(E_j,E_i)\\[6pt]

=-e^{2f}\big(\widetilde g(\widetilde R(E_1,E_i)E_j,E_1)+\widetilde
g(\widetilde R(E_1,E_j)E_i,E_1)\big)=0.
\end{array}
$$
Thus,
\begin{equation}\label{Ricci-i-j}
\widetilde\rho(E_i,E_j)=0 \quad {\rm{for}}~~  i,j\geq 2,~ i\neq j.
\end{equation}
Fix a point $p$ in the open set  $M'=\{q\in M: (\nabla f)(q)\neq
0\}$, and take an oriented $g$-orthonormal frame of vector fields
$E_1,E_2,...,E_n$ with $E_1=\frac{\nabla f}{||\nabla f||_g}$ in a
neighbourhood of $p$ in $M'$. Then $A_i=e^{-f}E_i$, $i=1,...,n$, is
an oriented $\widetilde g$-orthonormal frame. Let $X,Y,Z\in T_pM$
and extend $X,Y$ to vector fields in a neighbourhood of $p$ such
that $\widetilde\nabla X|_p=\widetilde\nabla Y|_p=0$. Set
$X_i=\widetilde g(X,A_i)$, $Y_i=\widetilde g(Y,A_i)$. Then, by
(\ref{Ricci 1-1/i}), (\ref{Ricci-i-i}) and (\ref{Ricci-i-j}),
$$
\begin{array}{c}
(\widetilde\nabla_{Z}\widetilde\rho)(X,Y)=\sum\limits_{i,j=1}^n
Z\big(X_{i}Y_{j}\,\widetilde\rho(A_i,A_j)\big)\\[6pt]

=\displaystyle{\frac{n-1}{2t}Z(X_1Y_1e^{-2f})+\frac{1}{2t}\sum\limits_{i=2}^n
Z(X_iY_ie^{-2f})}\\[6pt]

=\displaystyle{\frac{n-2}{2t}Z(X_1Y_1e^{-2f})+\frac{1}{2t}Z\big(\widetilde
g(X,Y)e^{-2f}\big)}\\[6pt]

=\displaystyle{\frac{n-2}{2t}Z(X_1Y_1e^{-2f})+\frac{1}{2t}\widetilde
g(X,Y)Z(e^{-2f})}.
\end{array}
$$
It follows that
$$
(\widetilde\delta\widetilde\rho)(Y)=-\sum_{i=1}^n
(\widetilde\nabla_{A_i}\widetilde\rho)(A_i,Y)=-\frac{1}{2t}[(n-2)A_1\big(g(A_i,Y)e^{-2f}\big)+Y(e^{-2f}].
$$
In particular,
$$
(\widetilde\delta\widetilde\rho)(A_1)=\frac{n-1}{t}e^{-2f}A_1(f).
$$
It follows from (\ref{Ricci 1-1/i}) and (\ref{Ricci-i-i}) that the
scalar curvature $\widetilde\tau$ of $\widetilde g$ is
$$
\widetilde\tau=\frac{n}{2t}e^{-2f}.
$$
Now, the well-know identity
$\widetilde\delta\widetilde\rho=-d\widetilde\tau$ (a consequence of
the differential Bianchi identity), implies
$$
\frac{n-1}{t}e^{-2f}A_1(f)=\frac{n}{t}e^{-2f}A_1(f).
$$
Therefore, $A_1(f)=0$ at the point $p$. At that point,
$A_1=e^{-f}\frac{\nabla f}{||\nabla f||_g}$. Hence $g(\nabla
f,\nabla f)=0$ at $p$. Thus, $(\nabla f)(p)=0$, a contradiction. It
follows that $f=const$ on $M$.

\smallskip

Conversely, if $f=const$, identity (\ref{Id harmonic}) is obvious,
hence the map $Id: ({\mathcal Z}, g_s)\to ({\mathcal Z},\widetilde
g_t)$ is harmonic.

\end{proof}

\smallskip

\noindent {\bf Acknowledgement}. I am very grateful to the reviewers
whose comments helped to correct and improve the final version of
the paper.

I would also like to thank Christian Yankov who has verified some of
the computations using  MAPLE software.

\end{document}